\documentclass[11pt,leqno]{article}
\usepackage{amsmath,amstext, amsthm,amsfonts,amssymb,amsthm} 
\numberwithin{equation}{section}
\def\C{\mathbb{ C}}

\def\D{\mathcal{D}}\def\L{\mathcal{L}}
\def\T{\mathcal{T}}\def\S{\mathcal{S}}\def\lam{\lambda}
\usepackage{epsfig}
\newtheorem{thm}{Theorem}[section]

\newtheorem{cor}[thm]{Corollary}

\newcommand{\be}{\begin{equation}}
\newcommand{\ee}{\end{equation}}
\newcommand\qbi[3]{{{#1}\atopwithdelims[]{#2}}_{#3}}
\newcommand{\ba}{\begin{array}}
\newcommand{\ea}{\end{array}}

\newcommand{\mbb}{\mathbb}
\newcommand{\mc}{\mathcal}

\newcommand{\gauss}[2]{\genfrac{[}{]}{0pt}{}{#1}{#2}_q}
\newcommand{\awp}{ Askey-Wilson polynomials }
\newtheorem{rem}[thm]{Remark}
\newcommand{\bea}{\begin{eqnarray}}
\newcommand{\eea}{\end{eqnarray}}

\newcommand{\Sum}{\sum_{n=0}^\infty}
\renewcommand{\(}{\left( } \renewcommand{\)}{\right) }

\title{Addition Theorems Via Continued Fractions}
    \author{Mourad E.H. Ismail
\\ Department of Mathematics \\ University of Central Florida
\\  Orlando, FL 32816\\
USA
\\ \and Jiang Zeng \\
Universit\'{e} de Lyon,  Universit\'{e} Lyon 1\\
Institute Camille Jordan,
UMR 5028 du CNRS\\
69622 Villeurbanne  \\
France}
\begin{document}
\maketitle
\begin{abstract}
We show connections between a special type of addition formulas and a theorem of Stieltjes and Rogers.  We use different techniques  to derive the desirable addition formulas. We apply our approach to derive special addition theorems for  Bessel functions and confluent hypergeometric functions. We also derive several additions theorems for basic hypergeometric functions.  Applications to the evaluation of Hankel determinants are also given .
\end{abstract}
\noindent {\it Mathematics Subject Classification} Primary 33D15, 33 C15,
Secondary 30E05, 05A15.
\noindent Keywords Addition theorems, orthogonal polynomials, continued $J$-fractions,
 $q$-orthogonal polynomials, Askey-Wilson polynomials,
 Bessel and confluent hypergeometric functions.
\maketitle
\section{Introduction}

An algebraic addition theorem for a function $f$ is an
identity of the form
\begin{equation}
\label{eqKoeKoo}
 P\bigl( f(x),\, f(y),\, f(x+y)\bigr) = 0
\end{equation}
for some polynomial $P$ in three variables.  Weierstrass  proved that
an analytic function satisfying an algebraic addition theorem
is a rational function in $z$, a rational function in
$e^{\lambda z}$ for some $\lambda$, or an elliptic
function,  \cite[Chapter 13]{For}. This notion is too restricted to be
useful in the theory of  special functions.
In general a family, say $\phi_\lambda$,
of special functions satisfies an addition formula
if there is an  elementary continuous function
$\Lambda$ of three variables $x,y,t$  and  an
expansion in terms of a  family
of special functions $\psi_\mu$ such that the expansion
coefficients factor as products in $x$ and $y$. In other word we have
\begin{equation}
\label{eq:generalformofadditionformula1}
\phi_\lambda\bigl(\Lambda(x,y,t)\bigr) =
\sum_\mu \, C(\lambda, \mu)\, \phi^\mu_\lambda(x) \,
\phi^\mu_\lambda(y)\, \psi_\mu(t), \qquad  C(\lambda, \mu)
\in \mathbb{C}.
\end{equation}

Recall the definition of a  Bessel function
\bea
\label{eqdefBesselJ}
J_{\nu}(z)=\sum_{m=0}^\infty\frac{(-1)^m}{m!\Gamma(m+\nu+1)}
\left(\frac{x}{2}\right)^{2m+\nu}=\frac{(z/2)^{\nu}}{\Gamma(\nu+1)}
{}_{0}F_{1}\left(\left.
\begin{array}{c}
-\\
\nu+1
\end{array}\right|\frac{-z^2}{4}
\right),
\eea
and the modified Bessel function
\bea
\label{eqdefBesselI}
I_\nu(x) := e^{-i\nu\pi/2}J_\nu (e^{i\pi/2}x) =
\sum_{m=0}^\infty\frac{(x/2)^{\nu+2m}}{m!\; \Gamma(m+\nu+1)}.
\eea
One important addition theorem is the addition theorem for Bessel functions,
\bea\label{eq:besseladd}
\frac{J_\nu(w)}{w^\nu} =  \frac{\Gamma(\nu)}{(zZ/2)^\nu}
\; \Sum (\nu+n) C_n^\nu(\cos \phi)\, J_{\nu+n}(z) \, J_{\nu+n}(Z),
\eea
for $\nu \ne 0,-1, -2, \cdots$, where $w := (z^2+Z^2-2zZ\cos \phi)^{1/2}$,
\cite[(7.15.30)]{Erd:Mag:Obe:Tri}, \cite{Wat}.  The polynomials $\{C_n^\nu(x)\}$ are the ultraspherical polynomials. The special case $\phi = \pi$ is
\bea
\label{eqBesselJ+}
\frac{J_\nu(x+y)}{(x+y)^\nu} =  \frac{\Gamma(\nu)}{(xy/2)^\nu}
\; \Sum (\nu+n) \frac{(-1)^n(2\nu)_n}{n!}\;  J_{\nu+n}(x) \, J_{\nu+n}(y),
\eea
since $C_n^\nu(-1) = (-1)^n(2\nu)_n/n!$. The addition theorems we will encounter in this work are of the type \eqref{eqBesselJ+}.

This work arose from an attempt to understand the Stieltjes-Rogers theorem of continued $J$-fractions, see Theorem \ref{St-Rog}. It is clear  that \eqref{eqaddthm} of Theorem \ref{St-Rog} is an addition theorem of the type  \eqref{eqBesselJ+}.

We decided to explore $q$-analogues of Rogers' addition formula \eqref{eqaddthm} of Theorem \ref{St-Rog} and  to compute the functions $Q_j(x), j =0, 1, \dots$ for specific
 continued fractions,  since the theory of orthogonal polynomials, especially the recently discovered one, provide a rich source of continued $J$-fractions.
We discovered two $q$-analogues of Theorem \ref{St-Rog}. They are Theorems~\ref{qaddthm} and \ref{NCaddthm}.

In this work we establish additions theorems for many special functions. To the best of our knowledge only \eqref{eqBesselJ+}  and \eqref{eq:addfor1F1}  are  known.  We offer three different techniques of proof and provide at least one example of each technique as  illustrations. We realize that it is possible to use fewer techniques to achieve the same goals but we believe there is merit in utilizing as many different ideas as possible.  One approach uses the plane wave expansion, \cite[(4.8.2)]{Ism2},
\begin{equation}
\label{eqPWJacobi}
\begin{split}
e^{xy} &=\sum_{n=0}^\infty\frac{\Gamma(\alpha+\beta+n+1)}{\Gamma(\alpha+\beta+2n+1)}\,
(2y)^ne^{-y}
\; {}_1F_1\(\left.\begin{matrix}
\beta+n+1 \\
\alpha+\beta+2n+2
\end{matrix}\,\right|2y\)\,P_n^{(\alpha,\beta)}(x),
\end{split}
\end{equation}
for $\alpha>-1$, $\beta>-1$,
and its special case \cite[(4.8.3)]{Ism2},
\begin{equation}
\label{eqPWUltra}
e^{xy}=\Gamma(\nu)(y/2)^{-\nu}
\sum_{n=0}^\infty (\nu+n)I_{\nu+n}(y)C_n^\nu(x), \quad \nu > -1/2.
\end{equation}
The polynomials $\{P_n^{(\alpha, \beta}(x)\}$ and $\{C_n^\nu(x)\}$ are Jacobi and ultraspherical polynomials, respectively. The expansions
\eqref{eqPWJacobi}--\eqref{eqPWUltra} are instances of the Fields and Wimp expansions \cite{Fie:Wim}, see also \cite{Fie:Ism}, \cite{Ver}.
Other techniques use Rodrigues type formulas followed by integration by parts, and connection coefficient formulas.

We shall follow the notation and terminology in \cite{And:Ask:Roy}, \cite{Ism2}, and \cite{Gas:Rah}. In particular we shall use the Rogers connection coefficients formula
for the continuous $q$-ultraspherical polynomials $\{C_n(x; \beta |q)\}$,
\begin{align}
\label{eq:connection}
C_n(x;\gamma|q)=\sum_{k=0}^{\lfloor {n/2}\rfloor}\frac{\beta^k(\gamma/\beta)_k(\gamma)_{n-k}}
{(q)_k(q\beta)_{n-k}}\frac{1-\beta q^{n-2k}}{1-\beta}C_{n-2k}(x;\beta|q),
\end{align}
\cite[p. 330]{Ism2}, and the facts that
\bea
\label{eqUnasultra}
U_n(x) = C_n^{1}(x) = C_n(x; q|q).
\eea
\section{Preliminaries}
Given a moment sequence $\{\mu_n\}$,
we define the linear functional ${\cal L}: x^n\mapsto \mu_n$ on the vector space of
polynomials $\C[x]$. We shall always assume $\mu_0 = {\cal L}(1) =1$. Then the monic polynomials $P_n(x)$ orthogonal with respect to the
$\cal L$ or the moment $\mu_n$ satisfy the following three term recurrence
relation (the  spectral theorem for orthogonal polynomials \cite[Chapter 2]{Ism2}):
\begin{align}\label{eq:threeterm}
P_{n+1}(x)=(x-b_n)P_{n}(x)-\lambda_n P_{n-1}(x),\quad n \ge 0,
\end{align}
where $\lambda_0 P_{-1}(x)=0$ and $P_0(x)=1$. We shall always require the functional to be regular, \cite{Chi},  which is equivalent to demanding that $\lambda_n \ne 0$ for all $n, n>0$.  The orthogonality relation is
\bea
\label{eqorthP}
{\cal L}(P_mP_n) = \lambda_1 \lambda_2 \cdots \lambda_n \delta_{m,n}.
\eea

The moment sequence is related to the coefficients $b_n$ and $\lambda_n$ by the following
identity:
\begin{equation}\label{eq:cf}
1+\sum_{n\geq 1}\mu _nx^n=
{1\over\displaystyle 1-b_0 x-
{\strut \lambda_1x^2\over\displaystyle 1-b_1x-
{\strut \lambda_2x^2\over\displaystyle
{\strut \ddots \over\displaystyle 1-b_nx-
{\strut \lambda_nx^2\over\displaystyle\ddots
}}}}}.
\end{equation}
Define the determinants
$$
\Delta_{i,n}=\left|
\begin{array}{cccc}
  \mu_0 & \mu_1 & \ldots & \mu_i \\
  \mu_1 & \mu_2 & \ldots & \mu_{i+1} \\
  \vdots & \vdots & \vdots & \vdots  \\
  \mu_{i-1} & \mu_{i} & \ldots & \mu_{2i-1} \\
  \mu_{n} & \mu_{n+1}  & \ldots & \mu_{n+i} \\
\end{array}
\right|,\qquad D_n(x)=\left|
\begin{array}{cccc}
  \mu_0 & \mu_1 & \ldots & \mu_n \\
  \mu_1 & \mu_2 & \ldots & \mu_{n+1} \\
  \vdots & \vdots & \vdots & \vdots  \\
  \mu_{n-1} & \mu_{n} & \ldots & \mu_{2n-1} \\
  1 & x  & \ldots & x^{n} \\
\end{array}
\right|.
$$
In particular, let
\bea\label{specialhankel}
D_n=\Delta_{n,n},\qquad \chi_n=\Delta_{n,n+1}.
\eea
Then $P_n(x)=(D_{n-1})^{-1}D_n(x)$ is the monic orthogonal polynomial sequence for $\cal L$.

It is easy to see that
\begin{align}
{\cal L}(x^nP_n(x))&=\frac{D_n}{D_{n-1}}=
\lambda_n\lambda_{n-1}\ldots \lambda_1,  \label{eqLxP}\\
{\cal L}(x^{n+1}P_n(x))&=\frac{\chi_n}{D_{n-1}}
=\lambda_n\lambda_{n-1}\ldots \lambda_1(b_0+\cdots +b_n).
\end{align}
Therefore
\begin{align}
\lambda_n=\frac{{\cal L}[P_n^2(x))]}{{\cal L}[P_{n-1}^2(x))]}
=\frac{D_{n-2}D_n}{D_{n-1}^2},
\end{align}
and
\begin{align}
b_n=\frac{{\cal L}[xP_n^2(x))]}{{\cal L}[P_n^2(x))]}
=\frac{\chi_n}{D_{n}}-\frac{\chi_{n-1}}{D_{n-1}}.
\end{align}
The next theorem is the backbone of this work.
\begin{thm}\label{St-Rog}
Define the Stieltjes tableau of entries $H_{i,n}$ $(i,n \ge 0)$ by
\begin{align}\label{stieltjes}
H_{i,n}&=0\quad \textrm{for}\quad i<0\quad \textrm{and}\quad i>n;\nonumber\\
H_{n,n}&=1\quad \textrm{for all }\quad n\geq 0;\\
H_{i,n}&=H_{i-1,n-1}+b_{i}H_{i,n-1}+\lambda_{i+1}H_{i+1,n-1}.\nonumber
\end{align}
Then the generating function $\sum_{n\geq 0}H_{0,n}x^n$
has the continued fraction expansion \eqref{eq:cf}
if and only if, for any two nonnegative integers $k,\ell\geq 0$,
 the following  convolution identities
\begin{align}\label{eq:stieltjes}
H_{0,k+\ell}=H_{0,k}H_{0,\ell}+\lambda_1H_{1,k}H_{1,\ell}+\cdots
+\lambda_1\cdots \lambda_j H_{j,k}H_{j,\ell}+\cdots,
\end{align}
hold. Moreover with the exponential generating functions of $\{Q_j(t)\}$
\begin{align}
Q_j(t)=\sum_{n= j}^\infty H_{j,n}\frac{t^n}{n!},
\end{align}
the convolution identity \eqref{eq:stieltjes} is equivalent  to the  addition formula
\bea
\label{eqaddthm}
Q_0(x+y)=\Sum \lambda_1\cdots \lambda_n\;  Q_n(x)Q_n(y).
\eea
\end{thm}
Wall \cite{Wal} points out that the first part of Theorem \ref{St-Rog} is due to Stieltjes but the addition theorem \eqref{eqaddthm} is due to Rogers. For a proof and references see
\cite[Section 53]{Wal}.

It is important to note  that the $H_{j,n}$'s are the  connection coefficients  in
\bea
\label{eqconnHjn}
x^n = \sum_{j=0}^n H_{j,n} P_j(x).
\eea
Observe  that the addition formula \eqref{eqaddthm} is equivalent to
\bea
\label{eqaddformh's}
\frac{xh_0(x)-yh_0(y)}{x-y} =\Sum \lambda_1\cdots \lambda_n\;  h_n(x)h_n(y),
\eea
where
$$
h_j(t)=\sum_{n= j}^\infty H_{j,n}t^n.
$$
In general \eqref{eq:stieltjes} implies
\bea
\sum_{m=0}^\infty H_{0,m} \sum_{j=0}^m c_js^jd_{m-j}t^{m-j} =
\Sum \lambda_1\cdots \lambda_n\;  Q_n(s)R_n(t),
\eea
where
\bea
Q_n(s) = \sum_{j=n}^\infty H_{n,j} c_js^j, \quad R_n(t) =\sum_{j=n}^\infty H_{n,j} d_jt^j.
\eea
We may define a generalized translation operator $(GT)$ on polynomials by
\bea
\label{eqdfnGT}
(GT)_s x^m = \sum_{j=0}^m c_js^jd_{m-j}t^{m-j}.
\eea
One can extend $(GT)_s$ to formal power series by linearity.

\begin{rem}
According to the Flajolet-Viennot theory \textup{\cite{Fla,Vie}} we can interpret $H_{i,n}$ in
the Stieltjes' tableau as follows.
Let us  attach weights to the steps of a lattice path at level $i$ $(i\geq 0)$
of a Motzkin path in the  following way:
$$
w(/)=1,\quad w(-)=b_{i}\quad \textrm{and}\quad w(\backslash)= \lambda_i.
$$
Let $\Gamma_{0\to i}(n)$ be the set of Motzkin paths from level 0 to level $i$ of length $n$.
Then we have the following interpretation:
\begin{align*}
H_{i,n}&=\sum_{\gamma\in \Gamma_{0\to i}(n)}w(\gamma),\\
\lambda_1\lambda_2\cdots \lambda_iH_{i,n}&=\sum_{\gamma\in \Gamma_{i\to 0}(n)}w(\gamma).
\end{align*}
This provides a combinatorial interpretation of \eqref{eq:stieltjes}.
\end{rem}

Let  $\{P_j(x)\}$  be  the monic  orthogonal polynomials with respect to the moment sequence
$H_{0,n}$ and ${\cal L}$ be the functional ${\cal L}(x^n)=H_{0,n}$.  Then it follows  from
\eqref{eqconnHjn} that
\begin{align}\label{eq:func}
Q_j(t)=\sum_{n = j}^\infty  \frac{{\cal L}\left(x^nP_j(x)\right)}
{{\cal L}\left(P_j^2(x)\right)}\frac{t^n}{n!}=\frac{1}{\lambda_1\ldots \lambda_j}{\cal L}\left(P_j(x)e^{xt}\right).
\end{align}
Note that $b_i=H_{i+1,i+2}-H_{i,i+1}$ and
\begin{align}\label{eq:hankel}
H_{j,n}=\frac{\Delta_{j,n}}{\lambda_1\ldots \lambda_j}.
\end{align}
Therefore the addition formula \eqref{eqaddthm} generalizes the Hankel
determinants in \eqref{specialhankel}. We refer the readers to
\cite{And:Wim,Hou:Las:Mu, Ism1, Kra1, Kra2,Wim} for the application of orthogonal polynomials
to the computation of Hankel determinants.

The following theorem gives another interpretation of the function $Q_0(x)$ for which our technique will derive an addition theorem.
\begin{thm}\label{thmcontour}
Assume that $\{P_n(x)\}$ are orthogonal with respect to a positive measure $\mu$ with compact support contained in $\{z: |z| <  r\}$. Then
\bea
\label{eqQj}
\quad \lambda_1 \lambda_2 \dots \lambda_j Q_j(t) = \oint_{|z|=r} e^{t/z}\, F_n(1/z)\, \frac{dz}{z},
\eea
 where
 \bea
\; F_n(z) := \int_{\mathbb {R}}\frac{P_n(u)}{z-u}d\mu(u), \quad z\notin \textup{supp}\{\mu\}.
\eea
\end{thm}
\begin{proof}  Let $\mu_n = \int_{\mathbb{R}} x^n d\mu(x)$.
The right-hand side of \eqref{eqQj} is
\bea
\begin{gathered}
 \oint_{|z|=r}\; e^{t/z}\int_{\mathbb {R}}\frac{P_j(u)}{1-zu} \, d\mu(u) \, \frac{dz}{z}
 = \int_{\mathbb {R}}   \oint_{|z|=r}  \Sum \frac{t^nz^{-n}}{n!} \sum_{k=0}^\infty(uz)^k P_j(u) \, d\mu(u) \, \frac{dz}{z}\\
 =\int_{\mathbb {R}} \Sum  \frac{(tu)^n }{n!}\, P_j(u) \, d\mu(u) =
 \int_{\mathbb {R}}  e^{tu} \, P_j(u) \, d\mu(u),
\end{gathered}
 \nonumber
\eea
 and the theorem follows from \eqref{eq:func}.
\end{proof}

The function $F_n(z)$ is related to the function of the second kind \cite{Ism2}.

Throughout this work we will use the Heine transformation
\bea
\label{eqHeine}
{}_2\phi_1(A, B, C, Z) = \frac{(B, AZ;q)_\infty}{(C, Z;q)_\infty}{}_2\phi_1(C/B, Z; AZ;q, B)
\eea
\cite[(III.1)]{Gas:Rah}
and the ${}_2\phi_1$ to ${}_2\phi_2$ transformation
\bea
\label{eq2phi2trans}
{}_2\phi_1(A, B;C;q,Z)= \frac{(AZ;q)_\infty}{(Z;q)_\infty}{}_2\phi_2(A, C/B; C, AZ; q, BZ).
\eea
\cite[(III.4)]{Gas:Rah}.

\section{Ultraspherical and Jacobi Polynomials}
The ultraspherical  (or Gegenbauer) polynomials are
$$
C_{n}^{\nu}(x)=\frac{(2\nu)_{n}}{n!}
{}_{2}F_{1}\left(\left.
\begin{array}{cccc}
-n,&n+2\nu\\
&\hspace{-1cm}\nu+ 1/2
\end{array}\right| \frac{1-x}{2}
\right), \quad \nu \neq 0.
$$
The normalized weight function is
$$
w(x)=(1-x^2)^{\nu-1/2}A(\nu),\quad A(\nu)=\frac{\Gamma(\nu+1)}{\Gamma(1/2)\Gamma(\nu+1/2)}.
$$
The corresponding orthogonality functional $\cal L$ is defined by
$$
{\cal L}(f)=\int_{-1}^1f(x)w(x)dx.
$$
The monic ultraspherical polynomials $\{P_n(x)\}$ and the $\lambda's$ are
$$
P_n(x)=\frac{n!}{(\nu)_n}\, 2^{-n}  \; C_n^{\nu}(x), \qquad
\lambda_j = \frac{j(j+2\nu-1)}{4(\nu+j-1)(\nu+j)}.
$$
Moreover
$$
xP_n(x)=P_{n+1}(x)+\frac{n(n+2\nu-1)}{4(\nu+n-1)(\nu+n)}P_{n-1}(x).
$$
Therefore \eqref{eq:func} when  $\nu > -1/2$ implies
\bea
\begin{gathered}
Q_i(x) = \frac{1}{\lambda_1 \lambda_2\dots \lambda_i}\;
 {\cal L}\left(e^{xt}P_i(x)\right)  \\
 = \frac{1}{\lambda_1 \lambda_2\dots \lambda_i} \Gamma(\nu)(t/2)^{-\nu}
 {\cal L}\left( \sum_{n=0}^\infty (\nu+n)I_{\nu+n}(t)C_n^\nu(x)P_i(x)\) \\
 = \Gamma(\nu)(t/2)^{-\nu}\, \frac{2^i \, (\nu)_i}{i!} (\nu+i) I_{\nu+i}(t).
\end{gathered}
\nonumber
\eea
Therefore
\bea
\label{eqQiforultra}
Q_i(t)=\frac{2^i\Gamma(\nu+i+1)}{i!(t/2)^\nu}I_{\nu+i}(t),
\eea
It is straightforward to see that addition formula \eqref{eqaddthm} in the present example
is equivalent to
\eqref{eqBesselJ+}.

Next we consider Jacobi polynomials.   The normalized weight function is
\begin{align}\label{eq:jacobiweight}
w(x)=(1-x)^{\alpha}(1+x)^\beta A(\alpha, \beta),\quad
A(\alpha, \beta)=\frac{\Gamma(\alpha+\beta+2)}
{2^{\alpha+\beta+1}\Gamma(\alpha+1)\Gamma(\beta+1)}.
\end{align}
In this case the functional to be considered  $\cal L$ is
$$
{\cal L}(f)=\int_{-1}^1f(x)w(x)dx.
$$
The monic Jacobi polynomials $\{P_n(x)\}$ are defined through
$$
P_n^{(\alpha,\beta)}(x)=\frac{(n+\alpha+\beta+1)_n}{2^nn!}P_n(x),
$$
so that
\begin{align}\label{eq:jacobiL}
\lambda_n=\frac{4n(n+\alpha)(n+\beta)(n+\alpha+\beta)}{(2n+\alpha+\beta-1)
(2n+\alpha+\beta)^2(2n+\alpha+\beta+1)}.
\end{align}
As in the case of ultraspherical polynomials we can use \eqref{eqPWJacobi}
and the result is
\bea
\label{eqQiJacobi}
Q_i(x)= \frac{t^ie^{-t}}{i!}{}_1F_1\left(\left. \begin{array}{c}
\beta+i+1\\
\alpha+\beta+2i+2
\end{array} \right|2t\right).
\eea
As an example of the use of Rodrigues formulas we give another derivation of
\eqref{eqQiJacobi}.  The Rodrigues formula for Jacobi polynomials is \cite[(4.2.8)]{Ism2}
\bea
(1-x)^\alpha (1+x)^\beta P_n^{(\alpha, \beta)}(x) =
\frac{(-1)^n}{2^n\; n!}  \frac{d^n}{dx^n} \left[(1-x)^{n+\alpha} (1+x)^{n+\beta} \right].
\eea
Therefore
\begin{align*}
Q_j(t)&=  \frac{2^j\, j! A(\alpha, \beta)}{\lambda_1\ldots \lambda_j (\alpha+\beta+j+1)_j}
\int_{-1}^1 e^{xt}\, (1-x)^{\alpha} (1+x)^{\beta} P_n^{(\alpha, \beta)}(x)\; dx
 \\
&= \frac{(-1)^j\, A(\alpha, \beta)}{\lambda_1\ldots \lambda_j (\alpha+\beta+j+1)_j}
\int_{-1}^1 e^{xt}\frac{d^j}{dx^j} \left[(1-x)^{j+\alpha} (1+x)^{j+\beta} \right] dx\\
&=  \frac{t^j\, A(\alpha, \beta)}{\lambda_1\ldots \lambda_j (\alpha+\beta+j+1)_j}
\int_{-1}^1 e^{xt}   (1-x)^{j+\alpha} (1+x)^{j+\beta}  dx,
\end{align*}
after integration by parts. Taking into account the integral representation
\cite[6.5.1)]{Erd:Mag:Obe:Tri-1}
\bea
{}_1F_1(a;c;z) = \frac{\Gamma(c)}{\Gamma(a)\Gamma(c-a)}
\; \int_0^1e^{zu} u^{a-1}(1-u)^{c-a-1}\, du, \;\; \;  \textup{Re} \; c > \textup{Re} \;  a >0,
\eea
\eqref{eq:jacobiL}, and \eqref{eq:jacobiweight},
we see that the last expression for  $Q_j(x)$ reduces to
\begin{equation}
\label{eqQiJacobi2}
Q_j(t)=\frac{t^j}{j!}\,e^t\,{}_1F_1\left(\left. \begin{array}{c}
\alpha+j+1\\
\alpha+\beta+2j+2
\end{array}\right| -2t\right).
\end{equation}

The equivalence of the representations \eqref{eqQiJacobi} and \eqref{eqQiJacobi2} follows from the transformation \cite[(1.4.11)]{Ism2}
$$
{}_1F_1\left(\left. \begin{array}{c}
a\\
c
\end{array}\right| z\right)
=e^z{}_1F_1\left(\left.
\begin{array}{c}
c-a\\
c
\end{array}\right| -z\right).
$$

This analysis establishes the following theorem, which is the main result of this section.
\begin{thm}\label{addition1F1}
We have the addition theorem for the confluent hypergeometric functions
\bea\label{eq:addfor1F1}
\begin{gathered}
{}_1F_1\left(\left. \begin{array}{c}
\alpha+1\\
\alpha+\beta+2
\end{array}\right| t+s \right) = \Sum  \frac{(\alpha+1)_n(\beta+1)_n(\alpha+\beta+1)_n}
{(\alpha+\beta+1)_{2n}(\alpha+\beta+2)_{2n}} \frac{(ts)^n}{n!}
\\
\qquad \qquad \times   {}_1F_1\left(\left. \begin{array}{c}
\alpha+n+1\\
\alpha+\beta+2n+2
\end{array}\right| t\right)\,{}_1F_1\left(\left. \begin{array}{c}
\alpha+n+1\\
\alpha+\beta+2n+2
\end{array}\right| s\right).
\end{gathered}
\eea
\end{thm}
When $\alpha = \beta= \nu-1/2$  Theorem \ref{addition1F1} reduces to \eqref{eqBesselJ+} since
\bea
e^{-x}\; {}_1F_1(\nu+1/2; 2\nu+1; 2x) = \Gamma(\nu+1)(2/x)^\nu I_\nu(x),
\eea
\cite[(6.9.10)]{Erd:Mag:Obe:Tri-1}.  Moreover  both
\eqref{eqQiJacobi} and \eqref{eqQiJacobi2} also reduce to
\eqref{eqQiforultra}.

Note that Theorem  \ref{addition1F1} and \eqref{eqQiforultra} can be proved from Theorem \ref{thmcontour} and the facts
\bea
\begin{gathered}
\int_{-1}^1 (1-t)^\alpha (1+t)^\beta \frac{P_n^{(\alpha, \beta)}(t)}{z-t}\, dt =
\frac{\Gamma(n+\alpha+1)\Gamma(n+\beta+1)}
{\Gamma(2n+\alpha+\beta+2)} \\
\qquad \times\frac{2^{n+\alpha+\beta+1}}{(z-1)^{n+1}}\;
 {}_2F_1\left(\left. \begin{array}{cc}
n+1, \alpha+n+1\\
\alpha+\beta+2n+2
\end{array}\right| \frac{2}{1-z}\right),
\end{gathered}
\nonumber
\eea
see (4.4.1) and (4.4.6) in \cite{Ism2}.  Moreover the integral
$$
\int_{-1}^1 e^{xt}\, (1-x)^\alpha(1+x)^\beta\, P_n^{(\alpha, \beta)}(x) dx
$$
can be evaluated from the plane wave expansion \eqref{eqPWJacobi} instead of the use of Rodrigues' formula.

It must be noted that \eqref{eq:addfor1F1}
coincides with formula (42) of Burchanl and Chaundy~\cite{Bur:Cha}. It is also
a limiting case of formula (50) in \cite{Bur:Cha}. The latter is
 stated in \cite{Erd:Mag:Obe:Tri-1},
see the first unnumbered formula after (7) in  Section   2.5.2.
Indeed if we replace $z$ by $z/b$ and $\zeta$ by
$\zeta/b$ and let $b \to \infty$, the above mentioned formula reduces to our
\eqref{eq:addfor1F1}.
The terminating case $\alpha=-m-1$, $\beta=\gamma+m$ of \eqref{eq:addfor1F1}
is the case $r=0$ of Koornwinder's addition formula for Laguerre polynomials,
see (3.3) in \cite{Koo}.
Also, this terminating case of \eqref{eq:addfor1F1} is the inverse of 10.12(42)
in \cite{Erd:Mag:Obe:Tri}.

\section{Two $q$-Addition Formulas}
The $q$-binomial formula
\begin{align}
\sum_{n=0}^\infty\frac{(a;q)_n}{(q;q)_n}t^n=\frac{(at;q)_\infty}{(t;q)_\infty},
\end{align}
yields Euler's $q$-analogues of exponential formula:
\begin{align}
\sum_{n=0}^\infty\frac{t^n}{(q;q)_n}&=\frac{1}{(t;q)_\infty}, \label{eqEuler1}\\
\sum_{n=0}^\infty\frac{q^{n\choose 2}t^n}{(q;q)_n}&=(-t;q)_\infty.\label{eqEuler2}
\end{align}
The $q$-difference operator is defined by
$$
\D_qf(z)=\frac{f(z)-f(qz)}{(1-q)z}.
$$
It is easy to see that
\begin{align}
\D_q^j((-xt;q)_\infty)&=\frac{t^jq^{j\choose 2}}{(1-q)^j}(-xtq^j;q)_\infty,\\
\D_q^j\left(\frac{1}{(xt;q)_\infty}\right)&=\frac{t^j(1-q)^{-j}}{(xt;q)_\infty}.
\end{align}
Note also that
\bea
\D_{1/q}^j\left(\frac{1}{(xt;q)_\infty}\right)
&=\frac{t^jq^{-{j\choose 2}}}{(1-q)^j}\frac{1}{(xtq^{-j};q)_\infty},\\
\D_{1/q}^j\left((-xt;q)_\infty\right)
&=\frac{t^j}{(1-q)^j}(-xt;q)_\infty.
\eea
\subsection{First $q$-Addition Formula}
Define the $q$-translation operator $\T_{s,q}$ by
\bea
\label{eqdfnofTy}
\T_{y,q} x^n = (x+y)(x+yq)\cdots (x+yq^{n-1}).
\eea
Extend this by linearity on functions $f(x)=\sum_{n= 0}^\infty  a_nx^n$, so that
$$
\T_{y,q}f(x)=\sum_{n = 0}^\infty a_n(x+y)(x+yq)\cdots (x+yq^{n-1}).
$$
Note that
\bea
\label{eqdfnofT/q}
\T_{y,q^{-1}}x^n=(x+y)(x+y/q)\cdots (x+y/q^{n-1})=\T_{x,q}y^nq^{-{n\choose 2}}.
\eea
Define two $q$-analogues of $Q_j(t)$ by
\bea
Q_j(t;q)=\sum_{n= j}^\infty H_{j,n}\frac{t^n}{(q;q)_n},\qquad
\widetilde{ Q}_j(t;q)=\sum_{n= j}^\infty  H_{j,n}\frac{q^{{n\choose 2}}t^n}{(q;q)_n}.
\eea

\begin{thm}\label{qaddthm}
The convolution identity \eqref{eq:stieltjes} is equivalent  to the
addition formula
\bea
\label{eqqaddthm}
\T_{s,q}Q_0(t;q)=\Sum \lambda_1\cdots \lambda_n\; Q_n(t;q)\;\widetilde{Q}_n(s;q).
\eea
Moreover, we have
\begin{align}
Q_j(t;q)&=\frac{1}{\lambda_1\cdots \lambda_j}{\cal L}
\left(\frac{P_j(x)}{(xt;q)_\infty}\right),\label{eq1}\\
\widetilde{Q}_j(t;q)&=\frac{1}{\lambda_1\cdots \lambda_j}
{\cal L}\left( P_j(x)(-xt;q)_\infty\right).\label{eq2}
\end{align}
\end{thm}
\begin{proof}
The $q$-binomial theorem and the definition of $\T_y$ gives
$$
\T_{y,q} x^n=(x+y)
(x+yq)\cdots (x+yq^{n-1})=\sum_{k=0}^n
{n\brack k}_qq^{k\choose 2}x^{n-k}y^k.
$$
This establishes the equivalence of \eqref{eqqaddthm} and \eqref{eq:stieltjes}.
Equations \eqref{eq1}  and \eqref{eq2} follow from Euler's formulas
 \eqref{eqEuler1}--\eqref{eqEuler2}.
\end{proof}

\subsection{Little Jacobi Polynomials} The little $q$-Jacobi
polynomials $\{p_n(x;a,b)\}$ are defined by
\bea
\label{eqdefnlittleqJ}
p_n(x;a,b) = {}_2\phi_1(q^{-n}, ab q^{n+1}; aq; q, qx).
\eea
 \cite[p.92-93]{Koe:Swa}. The corresponding monic polynomials
$\{P_j(x)\}$ are given by
\bea
p_n(x;a,b)=\frac{(-1)^nq^{-{n\choose 2}}(abq^{n+1};q)_n}{(aq;q)_n}P_n(x),
\eea
and
\bea
\lambda_n=\frac{aq^{2n-1}(1-q^n)(1-aq^n)(1-bq^n)(1-abq^n)}
{(1-abq^{2n-1})(1-abq^{2n})(1-abq^{2n})(1-abq^{2n+1})}.
\eea
Let
\bea
a=q^\alpha,\qquad b=q^\beta,
\eea
then the weight function is given by
$$
w(x;\alpha,\beta)=\frac{(qx;q)_\infty}{(q^{\beta+1}x;q)_\infty}x^\alpha,
$$
and the corresponding orthogonality functional is
\begin{align}
\label{eqfunctionaLQJ}
 {\cal L}(f)=\frac{(aq,bq;q)_\infty}{(abq^2,q;q)_\infty(1-q)}
\int_0^1f(x)w(x;\alpha,\beta)d_qx.
\end{align}
The Rodrigues-type formula for the little $q$-Jacobi polynomials is
\begin{align}\label{rodriguesqjacobi}
p_j(x;a,b)=\frac{1}{w(x;\alpha,\beta)}
\frac{(1-q)^jq^{j\alpha+{j\choose 2}}}
{(q^{\alpha+1};q)_j}  \D_{1/q}^j(w(x;\alpha+j,\beta+j)).
\end{align}
Therefore
\begin{align*}
{\widetilde Q}_j(t;q)&=\frac{a^{-j}q^{-j^2}(abq,abq^2;q)_{2j} (aq,bq;q)_\infty}
{(q,aq,bq,abq;q)_j (abq^2,q;q)_\infty (1-q)}
\int_0^1(-xt;q)_\infty w(x;\alpha,\beta) P_j(x)d_qx\\
&=\frac{a^{-j}q^{-j^2}(abq^2;q)_{2j} (aq,bq;q)_\infty}
{(q,bq;q)_j (abq^2,q;q)_\infty (1-q)}(-1)^j
q^{j\choose 2}\int_0^1(-xt;q)_\infty w(x;\alpha,\beta) p_j(x;a,b)d_qx\\
&=\frac{q^{-j}(abq^2;q)_{2j} (aq,bq;q)_\infty(-1)^j(1-q)^j}
{(q,aq,bq;q)_j (abq^2,q;q)_\infty (1-q)} \int_0^1(-xt;q)_\infty D_{1/q}^j(w(x;\alpha+j,\beta+j))d_qx.
\end{align*}
The $q$-analogue of integration by parts is
\begin{align}
\label{eqq-intgbyparts}
\int_a^b\D_q(f(x)) w(x)d_qx=-\frac{1}{q}\int_a^bf(x)\D_{1/q}(w(x))d_qx,
\end{align}
provided that  $w(a/q)=w(b/q)=0$, see \cite[(11.4.9)]{Ism2}. Applying
\eqref{eqq-intgbyparts} to the last expression for ${\widetilde Q}_j(t;q)$ we obtain
\begin{align*}
{\widetilde Q}_j(t;q)&=\frac{(abq^2;q)_{2j} (aq,bq;q)_\infty t^jq^{j\choose 2}}
{(q,aq,bq;q)_j  (abq^2,q;q)_\infty}
\sum_{n= 0}^\infty \frac{(-tq^{j+n},q^{n+1};q)_\infty}{(q^{\beta+n+1+j};q)_\infty}q^{(\alpha+j+1)n}\\
&=\frac{(-tq^j, aq^{j+1};q)_{\infty} t^jq^{j\choose 2}}
{(q;q)_j  (abq^{2+2j};q)_\infty}\;
{}_2\phi_1\left(\left.
\begin{array}{cc}
  bq^{j+1},&0   \\
   &\hspace{-1cm} -tq^j\\
\end{array}\right| q,aq^{j+1}
\right).
\end{align*}
This shows that
\bea
\label{eqQntilde}
{\widetilde Q}_j(t;q) =\frac{ t^jq^{j\choose 2}}
{(q;q)_j  }
\lim_{\epsilon \to 0}\frac{(-tq^j, aq^{j+1};q)_{\infty}}
{(abq^{2+2j},\epsilon;q)_\infty} \;
{}_2\phi_1\left(\left.
\begin{array}{cc}
  bq^{j+1},&\epsilon   \\
   &\hspace{-1cm} -tq^j\\
\end{array}\right| q,aq^{j+1}
\right).
\eea
Now the Heine transformation  \eqref{eqHeine}  leads to
$$
\frac{(-tq^j, aq^{j+1};q)_{\infty}}
{(abq^{2+2j},\epsilon;q)_\infty}
{}_2\phi_1\left(\left.
\begin{array}{cc}
  bq^{j+1},&\epsilon   \\
   &\hspace{-1cm} -tq^j\\
\end{array}\right| q,aq^{j+1}
\right)=
{}_2\phi_1\left(\left.
\begin{array}{cc}
  -tq^{j}/\epsilon,&aq^{j+1}   \\
   &\hspace{-1cm} abq^{2j+2}\\
\end{array}\right| q,\epsilon
\right).
$$
Therefore the above equation and \eqref{eqQntilde} establish the basic hypergeometric representation
\bea\label{tQjlittlejacobi}
{\widetilde Q}_j(t;q)=\frac{t^jq^{j\choose 2}}{(q;q)_j}\;
{}_1\phi_1\left(\left.
\begin{array}{c}
  aq^{j+1}  \\
  abq^{2j+2}\\
\end{array}\right| q, -tq^{j}
\right).
\eea
Similarly we have
\begin{align*}
Q_j(t;q)&=\frac{(abq^2;q)_{2j} (aq,bq;q)_\infty t^j}
{(q,aq,bq;q)_j  (abq^2,q;q)_\infty}
\sum_{n= 0}^\infty  \frac{(q^{1+n};q)_\infty a^nq^{(j+1)n}}{(bq^{n+1+j},tq^n;q)_\infty}\\
&=\frac{t^j (aq^{j+1};q)_\infty}
{(q;q)_j (t, abq^{2j+2};q)_{\infty}}\;
{}_2\phi_1\left(\left.
\begin{array}{cc}
  bq^{j+1},&t   \\
   &\hspace{-1cm} 0\\
\end{array}\right| q,aq^{j+1}
\right).
\end{align*}
Applying the Heine transformation \eqref{eqHeine}  to the last ${}_2\phi_1$  yields
\bea\label{Qjlittlejacobi}
Q_j(t;q)=\frac{t^j}{(q;q)_j}{}_2\phi_1\left(\left.
\begin{array}{cc}
  0,&aq^{j+1}   \\
   &\hspace{-1cm} abq^{2j+2}\\
\end{array}\right| q,t
\right).
\eea
\begin{thm}
For the functional in \eqref{eqfunctionaLQJ} the functions ${\widetilde Q}_j(t;q)$ and $Q_j(t;q)$ are given by  \eqref{tQjlittlejacobi} and \eqref{Qjlittlejacobi}. Moreover
we have the following addition formula:
\bea
\label{eqthm4.2}
\begin{gathered}
{}_2\phi_1\left(\left.
\begin{array}{cc}
  aq, & -s/t   \\
   &\hspace{-1cm} abq^{2}\\
\end{array}\right| q,t
\right) =\sum_{j=0}^\infty
\frac{q^{3\binom{j}{2}}(aq,bq,abq;q)_j}{(q;q)_j (abq,abq^2;q)_{2j}}\;  (qast)^j  \qquad \\
\qquad \qquad \times
{}_2\phi_1\left(\left.
\begin{array}{cc}
  0, &aq^{j+1}   \\
   &\hspace{-1cm} abq^{2j+2}\\
\end{array}\right| q,t
\right)\; {}_1\phi_1\left(\left.
\begin{array}{cc}
  aq^{j+1}   \\
   abq^{2j+2}\\
\end{array}\right| q, -s q^j
\right).
\end{gathered}
\eea
\end{thm}
\begin{proof}
We only need to show that the left-hand side of \eqref{eqthm4.2} is $\T_{s,q}Q_0(t,q)$. This is indeed the case as can be seen from \eqref{eqdfnofTy}.
\end{proof}

It is important to note that Theorem 4.2 is a $q$-analogue of Theorem 3.1.
Indeed with $a = q^\alpha, b = q^\beta$ and $s$ and $t$ replaced by $s(1-q)$ and $t(1-q)$, respectively, equation \eqref{eqthm4.2} tends to \eqref{eq:addfor1F1} as $q \to 1^-$.

Note that the transformation \eqref{eq2phi2trans} implies
\bea
Q_j(t;q) = \frac{t^j}{(q;q)_j(t;q)_\infty} {}_1\phi_1\left(\left.
\begin{array}{cc} bq^{j+1}   \\
   abq^{2j+2}\\
\end{array}\right| q, q^{j+1} at
\right).
\eea
Thus the addition theorem  \eqref{eqthm4.2} has the alternate form
 \bea
\label{eqthm4.2-2}
\begin{gathered}
{}_2\phi_1\left(\left.
\begin{array}{cc}
  aq, & -s/t   \\
   &\hspace{-1cm} abq^{2}\\
\end{array}\right| q,t
\right) =\sum_{j=0}^\infty
\frac{q^{3\binom{j}{2}}(aq,bq,abq;q)_j}{(q;q)_j (abq,abq^2;q)_{2j}}\;  \frac{(qast)^j}{(t;q)_\infty}  \quad \\
\qquad \qquad \times
{}_1\phi_1\left(\left.
\begin{array}{cc}
  bq^{j+1}   \\
  abq^{2j+2}\\
\end{array}\right| q, q^{j+1} at
\right)\; {}_1\phi_1\left(\left.
\begin{array}{cc}
  aq^{j+1}   \\
   abq^{2j+2}\\
\end{array}\right| q, - q^js
\right).
\end{gathered}
\eea
\noindent{\bf Remark}. A different $q$-analogue of Theorem 3.1  was given by
Jackson~\cite[(55)]{Jac}.

\subsection{Big $q$-Jacobi Polynomials}
The monic big $q$-Jacobi polynomials  are
\bea
p_n(x; a,b,c)=\frac{(abq^{n+1};q)_n}{(aq,cq;q)_n}P_n(x).
\eea
Let
\bea
w_1(x,a,b,c)=\frac{(x/a,x/c)_\infty}{(x,bx/c)_\infty},
\eea
and
\bea
w(x,a,b,c)=\frac{w_1(x,a,b,c)}{aq(1-q)}\frac{(aq,bq,cq,abq/c;q)_\infty}
{(q,abq^2,c/a,aq/c;q)_\infty}.
\eea
The corresponding orthogonality functional is
\bea
\label{eqBQJfunctional}
\L(f)=\int_{cq}^{aq}f(x)w(x,a,b,c)d_qx.
\eea
The Rodrigues-type formula for big $q$-Jacobi polynomials is
\bea
w_1(x)p_n(x; a,b,c)=\frac{a^nc^nq^{n(n+1)}(1-q)^n}{(aq,cq;q)_n}\D_q^n w_1(x,aq^n,bq^n,cq^n).
\eea
Note that
\bea
\lambda_n=\frac{-acq^{n+1}(1-q^n)(1-aq^n)(1-bq^n)(1-cq^n)(1-abq^n)(1-abq^n/c)}
{(1-abq^{2n-1})(1-abq^{2n})(1-abq^{2n})(1-abq^{2n+1})}.
\eea
So
\bea
\lambda_1\cdots \lambda_j=\frac{(-ac)^jq^{j(j+3)/2}
(q,aq,bq,cq, abq,abq/c;q)_j}
{(abq,abq^2;q)_{2j}}.
\eea
Therefore
\begin{align*}
Q_j(t, a,b,c)&=\frac{1}{\lam_1\cdots \lam_j}\L(P_j(x)/(xt;q)_\infty)\\
&=\frac{(-1)^j(1-q)^jq^{j+1\choose 2}(aq^{j+1},bq^{j+1},cq^{j+1},abq^{j+1}/c;q)_\infty}
{(q;q)_jaq(1-q)(q,abq^2,c/a,aq/c;q)_\infty}\\
&\qquad\times\int_{cq}^{aq}\frac{\D_q^jw_1(x,aq^j,bq^j,cq^j)}{(xt;q)_\infty}
d_qx\\
&=\frac{t^jq^{-j}(abq^2;q)_{2j}(aq^{j+1},bq^{j+1},cq^{j+1},abq^{j+1}/c;q)_\infty}
{a(q;q)_j(q,abq^2,c/a,aq/c;q)_\infty}I_j,
\end{align*}
where
\begin{align*}
I_j=&\frac{1}{q(1-q)}
\int_{cq}^{aq}\frac{(xq^{-j}/a,xq^{-j}/c;q)_\infty}{(x,
bx/c,xtq^{-j};q)_\infty}
d_qx\\
=&a\sum_{n= 0}^\infty  \frac{(q^{n+1}, aq^{n+1}/c;q)_\infty}
{(aq^{n+j+1}, abq^{n+j+1}/c;q)_\infty}\frac{q^{n+j}}{(atq^{n+1};q)_\infty}\\
&\qquad -c\sum_{n= 0}^\infty  \frac{(q^{n+1}, cq^{n+1}/a;q)_\infty}
{(cq^{n+j+1}, bq^{n+j+1};q)_\infty}\frac{q^{n+j}}{(ctq^{n+1};q)_\infty}\\
=&\frac{aq^j(aq/c,q;q)_\infty}{(aq^{j+1}, abq^{j+1}/c,
atq;q)_\infty} {}_3\phi_2\left(\left.
\begin{array}{ccc}
  aq^{j+1}, & abq^{j+1}/c, & atq \\
  aq/c, &  &\hspace{-2cm} 0 \\
\end{array}\right|q,q\right)\\
&-\frac{cq^j(cq/a,q;q)_\infty}{(cq^{j+1}, bq^{j+1}, ctq;q)_\infty}
{}_3\phi_2\left(\left.
\begin{array}{ccc}
  cq^{j+1}, & bq^{j+1}, & ctq \\
  cq/a, &  &\hspace{-2cm} 0 \\
\end{array}\right|q,q\right).
\end{align*}
The conclusion of the above calculations is that
\begin{align}\label{tQjbigjacobi}
\begin{gathered}
Q_j(t, a,b,c)
=\frac{t^j}{(q;q)_j(abq^{2j+2},q)_\infty}\hfill\\
\times \left[
\frac{(bq^{j+1}, cq^{j+1};q)_\infty}{( c/a, atq;q)_\infty}
{}_3\phi_2\left(\left.
\begin{array}{ccc}
  aq^{j+1}, & abq^{j+1}/c, & atq \\
  aq/c, &  &\hspace{-2cm} 0 \\
\end{array}\right|q,q\right)\right.\\
\left.+\frac{(aq^{j+1}, abq^{j+1}/c;q)_\infty}
{(a/c, ctq;q)_\infty}
{}_3\phi_2\left(\left.
\begin{array}{ccc}
  cq^{j+1}, & bq^{j+1}, & ctq \\
  cq/a, &  &\hspace{-2cm} 0 \\
\end{array}\right|q,q\right)\right].
\end{gathered}
\end{align}
Next we apply (12.5.8) in \cite{Ism2} with $A = D/(qat), B = aq^{j+1}, C = abq^{j+1}/c, E = abq^{2j+2}$ and let $D \to 0$ and  realize  that $Q_j(t, a,b,c)$ has the representation
\bea
\label{eqQjBigqJac}
\begin{gathered}
Q_j(t, a,b,c) = \frac{t^j}{(q;q)_j(aqt,q)_\infty}
{}_2\phi_1\left(\left.
\begin{array}{ccc}
  aq^{j+1},  abq^{j+1}/c\\
  abq^{2j+2} \\
\end{array}\right|q,qct \right).
\end{gathered}
\eea
It is clear from \eqref{eqQjBigqJac}  that $Q_j(t,a,b,c)=\frac{t^j}{(q;q)_j}+\cdots$.

Next we compute the functions ${\widetilde Q}_j(t,a,b,c)$. We have
\begin{align*}
{\widetilde Q}_j(t,a,b,c)&=
\frac{1}{\lam_1\cdots \lam_j}\L(P_j(x)(-xt;q)_\infty)\\
&=\frac{(-1)^jq^{j\choose 2}(1-q)^j(aq^{j+1},bq^{j+1},cq^{j+1},abq^{j+1}/c;q)_\infty}
{a(q;q)_j(q,abq^{2+2j},c/a,aq/c;q)_\infty}\\
&\times \int_{cq}^{aq}(-xt;q)_\infty\D_q^jw_1(x,aq^j,bq^j,cq^j)
d_qx\\
&=\frac{t^jq^{j\choose 2}(aq^{j+1},bq^{j+1},cq^{j+1},abq^{j+1}/c;q)_\infty}
{a(q;q)_j(q,abq^{2+2j},c/a,aq/c;q)_\infty}{\tilde I}_j,
\end{align*}
where
\begin{align*}
{\tilde I}_j&=\frac{q^{-j}}{q(1-q)}
\int_{cq}^{aq}\frac{(xq^{-j}/a,xq^{-j}/c;q)_\infty}{(x, bx/c,xtq^{-j};q)_\infty}
d_qx\\
&=aq^{-j}\sum_{n= 0}^\infty \frac{(-atq^{n+1},q^{n+1-j}, aq^{n+1-j}/c;q)_\infty}
{(aq^{n+1}, abq^{n+1}/c;q)_\infty}q^{n}\\
&\hspace{2cm}-cq^{-j}\sum_{n= 0}^\infty  \frac{(-ctq^{n+1}, q^{n+1-j}, cq^{n+1-j}/a;q)_\infty}
{(cq^{n+1}, bq^{n+1};q)_\infty} q^{n}\\
&=\frac{a(-atq^{j+1},q, aq/c;q)_\infty}
{(aq^{j+1}, abq^{j+1}/c;q)_\infty}
{}_3\phi_2\left(\left.
\begin{array}{ccc}
  aq^{j+1}, & abq^{j+1}/c, & 0 \\
  -atq^{j+1} &  &\hspace{-2cm} aq/c\\
\end{array}\right|q,q\right)\\
&-\frac{c(-ctq^{j+1},q, cq/a;q)_\infty}{(cq^{j+1}, bq^{j+1};q)_\infty}
{}_3\phi_2\left(\left.
\begin{array}{ccc}
  cq^{j+1}, & bq^{j+1}, & 0 \\
  -ctq^{j+1}, &  &\hspace{-2cm} cq/a\\
\end{array}\right|q,q\right).
\end{align*}
It follows that
\begin{align}\label{eqQjbigjacobi}
\begin{gathered}
{\widetilde Q}_j(t,a,b,c)
=\frac{t^jq^{j\choose 2}}{(q;q)_j(abq^{2+2j};q)_\infty}\hfill\\
\times \left[\frac{(bq^{j+1},cq^{j+1}, -atq^{j+1};q)_\infty}
{(c/a;q)_\infty}
{}_3\phi_2\left(\left.
\begin{array}{ccc}
  aq^{j+1}, & abq^{j+1}/c, & 0 \\
  -atq^{j+1} &  &\hspace{-2cm} aq/c\\
\end{array}\right|q,q\right)
\right.\\
\left.+\frac{(aq^{j+1},abq^{j+1}/c, -ctq^{j+1};q)_\infty}
{(a/c;q)_\infty}
{}_3\phi_2\left(\left.
\begin{array}{ccc}
  cq^{j+1}, & bq^{j+1}, & 0 \\
  -ctq^{j+1} &  &\hspace{-2cm} cq/a
\end{array}\right|q,q\right)\right].
\end{gathered}
\end{align}
To simplify  \eqref{eqQjbigjacobi} we apply (12.5.8) in \cite{Ism2} with
$B = aq^{j+1}, C = abq^{j+1}/c, D = -atq^{j+1}, E = abq^{2j+2}$ and let $A \to \infty$. Alternately we may use the result in Exercise 3.8 of \cite{Gas:Rah}. The conclusion is that
\bea
\label{eqQjbigjacobi2}
\begin{gathered}
{\widetilde Q}_j(t,a,b,c) = \frac{q^{\binom{j}{2}}t^j(-atq^{j+1};q)_\infty}{(q;q)_j}
 {}_2\phi_2\left(\left.
\begin{array}{ccc}
  aq^{j+1}, & abq^{j+1}/c \\
abq^{2j+2}, &   -atq^{j+1} \\
\end{array}\right|q,-tcq^{j+1}\right).
\end{gathered}
\eea
Finally the transformation \eqref{eq2phi2trans} gives yet another alternate representation, namely
\bea
\label{eqQjbigjacobi3}
{\widetilde Q}_j(t,a,b,c) =  \frac{q^{\binom{j}{2}}t^j}{(q;q)_j}\; (-t;q)_\infty\;
 {}_2\phi_1\left(\left.
\begin{array}{ccc}
  aq^{j+1},  cq^{j+1} \\
abq^{2j+2}  \\
\end{array}\right|q,-t \right).
\eea

It is clear from \eqref{eqQjbigjacobi2} or \eqref{eqQjbigjacobi3} that ${\widetilde Q}_j(t,a,b,c)=\frac{t^jq^{j\choose 2}}{(q;q)_j}+\cdots$.
\begin{thm}   The functions $Q_j(t,a,b,c)$ and ${\widetilde Q}_j(t,a,b,c)$ associated with
 the big $q$-Jacobi functional \eqref{eqBQJfunctional}  are
defined in \eqref{eqQjBigqJac}  and \eqref{eqQjbigjacobi3} \textup{(}or in  \eqref{eqQjbigjacobi2}\textup{)}. Moreover
we have the following addition formula:
\bea
\label{eqaddbigqJacobi}
\begin{gathered}
\frac{(-qas;q)_\infty}{(- s;q)_\infty}\;  {}_3\phi_2\left(\left.
\begin{array}{ccc}
  qa, & q ab/c & -s/t \\
abq^2, &   -qas \\
\end{array}\right|q, qct\right)  \\
=\sum_{j=0}^\infty
\frac{(-acst)^j q^{j(j+1)}
(aq,bq,cq, abq,abq/c;q)_j}
{(q;q)_j(abq,abq^2;q)_{2j}}\\
\times {}_2\phi_1\left(\left.
\begin{array}{ccc}
  aq^{j+1},  abq^{j+1}/c\\
  abq^{2j+2} \\
\end{array}\right|q,qct \right)\;  {}_2\phi_1\left(\left.
\begin{array}{ccc}
  aq^{j+1},  cq^{j+1} \\
abq^{2j+2}  \\
\end{array}\right|q,-s \right).
\end{gathered}
\eea
\end{thm}
\begin{proof}
We only need to evaluate $\T_{s,q} Q_0(t;a,b,c)$. Clearly $\T_{s,q} Q_0(t;a,b,c)$ is
\begin{align*}
&\Sum \frac{(qa, qab/c;q)_n}{(q, abq^2;q)_n}(qc)^n
 \sum_{k=0}^\infty \frac{(qa)^k}{(q;q)_k} \T_{s,q} t^{n+k} \\
 & =  \Sum \frac{(qa, qab/c;q)_n}{(q, abq^2;q)_n}(qc)^n
 \sum_{k=0}^\infty \frac{(qa)^k}{(q;q)_k} t^{n+k}(-s/t;q)_{n+k} \\
 &= \Sum \frac{(qa, qab/c-s/t ;q)_n}{(q, abq^2;q)_n}(qct)^n \; {}_1\phi_0(-q^ns/t;--;q, qat)\\
 &= \Sum \frac{(qa, qab/c-s/t ;q)_n}{(q, abq^2;q)_n}(qct)^n \frac{(-q^{n+1}as;q)_\infty}{(aqt;q)_\infty}.
\end{align*}
 The $q$-binomial theorem reduces the  last expression to
 \bea
\frac{(-qas;q)_\infty}{(qat;q)_\infty} \; {}_3\phi_2\left(\left.
\begin{array}{ccc}
  qa, & q ab/c & -s/t \\
abq^2, &   -qas \\
\end{array}\right|q, qct\right),
\nonumber
 \eea
 and the theorem follows.
\end{proof}

\subsection{Second $q$-Addition Formula}
Define the non-commutative operation $\S_y$ on $f(x)=\sum_{n= 0}^\infty a_nx^n$ by
$$
\S_yf(x)=\sum_{n = 0}^\infty  a_n(x+y)^n,
$$
where $yx=qxy$. Recall~\cite[p.28]{Gas:Rah} that the non-commutative binomial theorem reads
\begin{align}\label{eq:NCbinomial}
(x+y)^n=\sum_{k=0}^n{n\brack k}_qx^ky^{n-k}.
\end{align}

\begin{thm}\label{NCaddthm}
The convolution identity \eqref{eq:stieltjes} is equivalent  to the
non-commutative addition formula
\bea
\label{eqNCaddthm}
\S_sQ_0(t;q)=\Sum \lambda_1\cdots \lambda_n\; Q_n(t;q)\;Q_n(s;q).
\eea
Moreover, we have
\begin{align}
Q_j(t;q)&=\frac{1}{\lambda_1\cdots \lambda_j}{\cal L}
\left(\frac{P_j(x)}{(xt;q)_\infty}\right).\label{eq3}
\end{align}
\end{thm}
\begin{proof}
The non-commutative binomial theorem and the definition of $S_y$ gives
$$
S_y x^n=\sum_{k=0}^n
{n\brack k}_qx^{n-k}y^k.
$$
Multiply \eqref{eq:stieltjes} by $t^ms^n/(q;q)_m(q;q)_n$ and sum over all $m,n\geq 0$ we get
\begin{align*}
\sum_{m,n\geq 0}H_{0,m+n}\frac{t^ms^n}{(q;q)_m(q;q)_n}&=
\sum_{j,m,n\geq 0}\frac{t^ms^n}{(q;q)_m(q;q)_n}H_{j,m}H_{j,n}\lambda_1\ldots
\lambda_j.
\end{align*}

This establishes the equivalence of \eqref{eqNCaddthm} and \eqref{eq:stieltjes}.
Equations \eqref{eq3} follows from Euler's formula.
\end{proof}

\noindent{\bf Example}: We consider the  Rogers-Szeg\H{o} polynomials, which
are defined by
\begin{align}\label{eq:rogers-szego-poly}
h_n(a;q)=\sum_{k=0}^n{n\brack k}_q a^k.
\end{align}
They have the $q$-exponential generating function
$$
Q_0(t;q)=\sum_{n= 0}^\infty  \frac{h_n(a;q)}{(q;q)_n}t^n=\frac{1}{(t;q)_\infty (at;q)_\infty}.
$$
Therefore
\begin{align*}
\S_sQ_0(t;q)&=\sum_{n=0}^\infty \frac{h_n(a;q)}{(q;q)_n}(t+s)^n\\
&=\frac{1}{(t+s;q)_\infty (a(t+s);q)_\infty}.
\end{align*}
The corresponding orthogonal polynomials are Al-Salam-Carlitz polynomials $U_n^{(a)}(x;q)$, which have the generating function
\begin{align}
\label{eq:alsalamcarlitzgf}
 \frac{(t,at;q)_\infty}{(xt;q)_\infty}=
\sum_{n=0}^\infty\frac{U_n^{(a)}(x;q)}{(q;q)_n}t^n.
\end{align}
The associated functional is defined by
\begin{align}
\label{eqfunctionalAl-Ca}
 {\cal L}f(x)=\frac{1}{(1-q)(q,a,q/a;a)_\infty}\int_a^1(qx,qx/a;q)_\infty f(x)d_qx.
\end{align}
Note that \cite{Koe:Swa}
\bea
\label{eqlambdaAl-Ca}
\lambda_j=-aq^{j-1}(1-q^j).
\eea
It follows from \eqref{eq:alsalamcarlitzgf} that the generating function of $Q_n(x)$ is
\begin{align*}
\sum_{n= 0}^\infty \; \lambda_1\cdots \lambda_nQ_n(t;q)\frac{y^n}{(q;q)_n}&=
\frac{(y,ay;q)_\infty}{(1-q)(q,a,q/a;q)_\infty}\int_a^1
\frac{(qx,qx/a;q)_\infty}{(xy,xt;q)_\infty}d_qx\\
&=
\frac{(ay;q)_\infty}{(a,t;q)_\infty}\;
{}_2\phi_1\left(\left.
\begin{array}{cc}
  y,&t   \\
   &\hspace{-1cm} q/a\\
\end{array}\right| q,q
\right)+
\frac{(y;q)_\infty}{(1/a,at;q)_\infty}\;
{}_2\phi_1\left(\left.
\begin{array}{cc}
  ay,&at   \\
   &\hspace{-1cm} aq\\
\end{array}\right| q,q
\right)\\
&=\frac{(ayt,yq/t,q/ayt;q)_\infty}{(t,at,q/t,q/at;q)_\infty}\;
{}_2\phi_1\left(\left.
\begin{array}{cc}
  y,&ay   \\
   &\hspace{-1cm} yq/t\\
\end{array}\right| q,q/ayt
\right),
\end{align*}
where the last equality follows from
the transformation~\cite[(III.32)]{Gas:Rah}) if we replace all the small letters by capital ones and apply the parameter identification:
$$
A=y,\quad B=ay,\quad C=yq/t,\quad Z=q/ayt.
$$
Therefore
\begin{align*}
 &\sum_{n= 0}^\infty    \lambda_1\cdots \lambda_nQ_n(t;q)\frac{y^n}{(q;q)_n} =
\frac{(ayt,yq/t,q/ayt;q)_\infty}{(t,at,q/t,q/at;q)_\infty}\lim_{\delta \to 1^-}\;
{}_2\phi_1\left(\left.
\begin{array}{cc}
  y,&ay   \\
   &\hspace{-1cm} yq/t\\
\end{array}\right|q, q\delta/ayt
\right)\\
&=\frac{(ayt,yq/t;q)_\infty}{(t,at,q/t,q/at;q)_\infty}\lim_{\delta \to 1^-}(\delta q;q)_\infty\; {}_2\phi_1\left(\left.
\begin{array}{cc}
  q/t,&q/at   \\
   &\hspace{-1cm} yq/t\\
\end{array}\right|q, \delta
\right)\\
&=\frac{(ayt;q)_\infty}{(t,at;q)_\infty},
\end{align*}
where we used the transformation \cite[(III.3)]{Gas:Rah} and $\lim_{\delta\to 1^-}(1-\delta)\sum_{n=0}^\infty a_n \delta^n=\lim_{n\to \infty}a_n$.
By equating the coefficients of $y^n$ we get
\begin{align}
Q_n(t;q)=\frac{t^n}{(q;q)_n}\frac{1}{(t,at;q)_\infty}.
\end{align}
Summarizing we get the following addition formula.
\begin{thm} If $st=qts$ then the following addition formula holds
\begin{align}
\frac{1}{(t+s, a(t+s);q)_\infty}=\sum_{j=0}^\infty (-a)^jq^{j\choose 2}(q;q)_j
\frac{t^j}{(q;q)_j(t,at;q)_\infty}\frac{s^j}{(q;q)_j(s,as;q)_\infty}.
\end{align}
\end{thm}

\subsection{Computing  First $q$-Addition Formulas Using Generating Functions}
Recall that the  Al-Salam-Carlitz polynomials
$\{U_n^{(a)}(x;q)\}$   have the generating function \eqref{eq:alsalamcarlitzgf}
and the associated functional is in \eqref{eqfunctionalAl-Ca}.
It follows from \eqref{eq:alsalamcarlitzgf} that the generating function of $Q_n(t;q)$ is
\begin{align*}
\sum_{n= 0} ^\infty \lambda_1\cdots \lambda_nQ_n(t;q)\frac{y^n}{(q;q)_n}&=
\frac{(y,ay;q)_\infty}{(1-q)(q,a,q/a;q)_\infty}\int_a^1
\frac{(qx,qx/a;q)_\infty}{(xy,xt;q)_\infty}d_qx\\
&=
\frac{(ay;q)_\infty}{(a,t;q)_\infty}\;
{}_2\phi_1\left(\left.
\begin{array}{cc}
  y,&t   \\
   &\hspace{-1cm} q/a\\
\end{array}\right| q,q
\right)+
\frac{(y;q)_\infty}{(1/a,at;q)_\infty}\;
{}_2\phi_1\left(\left.
\begin{array}{cc}
  ay,&at   \\
   &\hspace{-1cm} aq\\
\end{array}\right| q,q
\right)\\
&=\frac{(ayt,yq/t,q/ayt;q)_\infty}{(t,at,q/t,q/at;q)_\infty} \;
{}_2\phi_1\left(\left.
\begin{array}{cc}
  y,&ay   \\
   &\hspace{-1cm} yq/t\\
\end{array}\right| q,q/ayt
\right),
\end{align*}
where the last equality follows from
the transformation~\cite[(III.32)]{Gas:Rah}) if we replace all the small letters by capital ones and then take the following substitutions:
$$
A=y,\quad B=ay,\quad C=yq/t,\quad Z=q/ayt.
$$
Therefore
\begin{align*}
 \sum_{n= 0}^\infty \lambda_1\cdots \lambda_nQ_n(t;q)\frac{y^n}{(q;q)_n} &=
\frac{(ayt,yq/t,q/ayt;q)_\infty}{(t,at,q/t,q/at;q)_\infty}\lim_{\delta \to 1} \;
{}_2\phi_1\left(\left.
\begin{array}{cc}
  y,&ay   \\
   &\hspace{-1cm} yq/t\\
\end{array}\right|q, q\delta/ayt
\right)\\
&=\frac{(ayt;q)_\infty}{(t,at;q)_\infty}.
\end{align*}
Equating the coefficients of $y^n$ in the above identity we get
\begin{align}
Q_n(t;q)=\frac{t^n}{(q;q)_n}\frac{1}{(t,at;q)_\infty}.
\end{align}
Similarly we have
\begin{align}
\label{eqGFofQtilde}
\sum_{n= 0}^\infty \lambda_1\cdots \lambda_n\widetilde Q_n(t;q)\frac{y^n}{(q;q)_n}&=
\frac{(y,ay;q)_\infty}{(1-q)(q,a,q/a;q)_\infty}
\int_a^1\frac{(qx,qx/a;q)_\infty}{(xy;q)_\infty}(-xt;q)_\infty d_qx.
\end{align}
The right-hand side of \eqref{eqGFofQtilde} is
\bea
\begin{gathered}
\frac{(y,ay;q)_\infty}{(a,q/a;q)_\infty}
\left[\Sum \frac{( q^{n+1}/a, -tq^n;q)_\infty}
{(q;q)_n\, (yq^n;q)_\infty)} q^n -  a \Sum \frac{(aq^{n+1},  -taq^n;q)_\infty}
{(q;q)_n\, (ayq^n;q)_\infty} q^n
 \right]\\
 = \frac{(ay,-t;q)_\infty}{(a;q)_\infty} \; {}_3\phi_2\left(\left.
\begin{array}{cc}
  y, 0, 0   \\
 -t, q/a\\
\end{array}\right|q, q
\right) +    \frac{(y, -at;q)_\infty}{(1/a;q)_\infty} \; {}_3\phi_2\left(\left.
\begin{array}{cc}
  ay, 0, 0   \\
 qa, -a t \\
\end{array}\right|q, q
\right).
\end{gathered}
\nonumber
\eea
In  (III.34) of \cite{Gas:Rah} replace $a,b,c,d,e$ by $A,B,C, D, EC$, respectively then let $C\to 0$ then let $A\to \infty$. The result is the three term relation
\bea
\begin{gathered}
{}_1\phi_1\left(\left.
\begin{array}{cc}
  B   \\
 D\\
\end{array}\right|q, \frac{DE}{B}
\right)= \frac{(E;q)_\infty}{(E/B;q)_\infty}
\; {}_3\phi_2\left(\left.
\begin{array}{cc}
  B, 0, 0   \\
 D, qB/E\\
\end{array}\right|q, q
\right)  \qquad \qquad  \\
\qquad \qquad + \frac{(B, DE/B;q)_\infty}{(D, B/E;q)_\infty}
\; {}_3\phi_2\left(\left.
\begin{array}{cc}
  E,  0, 0   \\
 DE/B, qE/B\\
\end{array}\right|q, q  \right).
\end{gathered}
 \nonumber
\eea
We now choose $B = y, E = ay, D =-t$. Therefore \eqref{eqGFofQtilde} becomes
\bea
\sum_{n= 0}^\infty \lambda_1\cdots \lambda_n\widetilde Q_n(t;q)\frac{y^n}{(q;q)_n}=
(-t;q)_\infty \; {}_1\phi_1\left(\left.
\begin{array}{cc}
  y   \\
 -t \\
\end{array}\right|q, -at
\right)
\eea
Finally we apply the $q$-binomial theorem in the form
$$
(y;q)_n = \sum_{k=0}^n \frac{(q;q)_n}{(q;q)_k(q;q)_{n-k}} q^{\binom{k}{2}} (-y)^k
$$
and \eqref{eqlambdaAl-Ca} to obtain
\bea
\widetilde Q_j(t;q) = \frac{(-tq^j;q)_\infty}{(q;q)_j} \; t^jq^{\binom{j}{2}} \;
{}_1\phi_1(0; -tq^j;q, -atq^j).
\eea
\begin{thm}
We have the identity
\bea
\label{eqAl-Ca-identity}
\begin{gathered}
\frac{(-s;q)_\infty}{(t;q)_\infty}\; {}_2\phi_1\left(\left.
\begin{array}{cc}
  0,-s/t   \\
  -s\\
\end{array}\right|q, at
\right)    \qquad \qquad \qquad \qquad \\
\qquad = \frac{1}{(t,at;q)_\infty} \Sum \frac{(-sq^n;q)_\infty}{(q;q)_n}
(-ast)^nq^{n(n-1)} \; {}_1\phi_1(0; -tq^n;q, -atq^n).
\end{gathered}
\eea
\end{thm}
\begin{proof}We need only to show that the left-hand side in \eqref{eqAl-Ca-identity} is
$\T_{s,q}Q_0(t;q) $. This can be seen as follows. By \eqref{eq:rogers-szego-poly} and
\eqref{eqdfnofTy} we have
\begin{align*}
\T_{s,q}Q_0(t;q)
&=\sum_{n=0}^\infty \sum_{k=0}^n
\frac{t^na^k(-s/t;q)_n}{(q;q)_k(q;q)_{n-k}} \\
&=\sum_{k= 0}^\infty \frac{t^ka^k(-s/t;q)_k}{(q;q)_k}\sum_{n= 0}^\infty
\frac{t^n(-sq^k/t;q)_n}{(q;q)_n}\\
&=\frac{(-s;q)_\infty}{(t;q)_\infty}\; {}_2\phi_1\left(\left.
\begin{array}{cc}
  0,&-s/t   \\
   &\hspace{-1cm} -s\\
\end{array}\right|q, at
\right).
\end{align*}
\end{proof}

\section{Sheffer-type Polynomials}
\subsection{Moments of Sheffer-type Polynomials}
For the Sheffer type orthogonal polynomials
 we can compute $Q_j(t)$ by the generating function
\begin{align}\label{eq:funcgen}
\Sum \lambda_1\ldots \lambda_n \alpha_n Q_n(t)\frac{y^n}{n!}=
{\cal L}\left(e^{xt}\Sum \alpha_nP_n(x)\frac{y^n}{n!}\right),
\end{align}
where $\alpha_j$'s are some suitably chosen constants.

The cases of  Hermite Laguerre,  Meixner, and Charlier polynomials do not
lead to interesting addition theorems because the addition theorem predicted
by Theorem \ref{St-Rog} follow from the binomial theorem. We just indicate the
corresponding $Q_j(t)$ in the addition formula of each family of polynomials.

Indeed in the case of
 Hermite polynomials $\{H_n(x)\}$,
\bea
\lambda_n = n/2, \qquad H_n(x) = 2^n P_n(x),
\eea
and
\bea
\label{eqQiforHer}
Q_n(t)= \frac{t^n}{n!}\; \exp(t^2/4).
\eea
In the case of Laguerre polynomials $\{L_n^\alpha(x)\}$,
\bea\label{laguerre}
\lambda_n = n(\alpha+n), \quad L_n^\alpha(x) = \frac{(-1)^n}{n!} P_n(x),
\eea
and
\bea
\label{eqQiforLag}
Q_n(t)=   \frac{t^n}{n!}\; (1-t)^{-\alpha-n-1}.
\eea
For Meixner polynomials $\{M_n(x;\beta,c)\}$, we have
\bea
\lambda_n = \frac{n (n+\beta-1)c}{(1-c)^2}, \quad
M_n(x;\beta,c) = \frac{(c-1)^n}{c^n\, (\beta)_n}\; P_n(x),
\eea
and
\bea
\label{eqQiforMeix}
Q_n(t)=\left(\frac{1-c}{1-ce^t}\right)^{\beta+n}
\frac{(e^t-1)^n}{n!}.
\eea
In the case of Charlier polynomials $\{C_n(x;a)\}$
\bea
\lambda_n = an, \qquad  C_n(x;a) = (-a)^{-n} P_n(x).
\eea
A calculation gives
\bea
\label{eqQiforChar}
Q_j(x) = \frac{(e^t-1)^j} {j!} \; \exp(e^t-1).
\eea
In the case of Meixner-Pollaczek polynomials $P_n^{\lambda}(x;\phi)$, \cite{Koe:Swa},
\bea
\lambda_n=\frac{n(n+2\lambda-1)}{4\sin^2\phi},\qquad
P_n^{\lambda}(x;\phi)=\frac{(2\sin \phi)^n}{n!}P_n(x).
\eea
One can see that the $Q_j$'s are given by
\bea
\label{eqQiforMP}
Q_j(x) = \frac{2^j}{j!}
 \; \left(\frac{\sin\phi}{\sin(t/2+\phi))}\right)^{2\lambda+j}\; \left[\sin (t/2)\right]^j.
\eea
Note that for the orthogonal  polynomials of Sheffer type
all the $Q_0(t)$'s have been given in \cite{Zen}.

\subsection{Sheffer-type Polynomials as Moments}
For the Hermite polynomials,  we have
\begin{align}
\exp{(2xt-t^2)}=\sum_{n= 0}^\infty \frac{H_{n}(x)}{n!}t^n.
\end{align}
Let
$$Q_{n}(t)=\frac{t^n}{n!}e^{2xt-t^2}=\frac{t^n}{n!}
+2(n+1)x\frac{t^{n+1}}{(n+1)!}+\cdots.
$$
From $\exp{(2x(t+s)-(t+s)^2)}=\exp{(2xt-t^2+2xs-s^2-2ts)}$
we derive the addition formula.
\begin{thm}
The functions $\{Q_j(t) \}$ have the addition formula
\begin{align}\label{addhermite}
Q_0(t+s)=\sum_{n= 0}^\infty n!(-2)^nQ_{n}(t)Q_{n}(s).
\end{align}
\end{thm}
It follows that $\lambda_n=-2n$ and
\begin{align}
H_{i,i+n}={i+n\choose i}H_n(x).
\end{align}
{\bf Remark}. Radoux~\cite{Rad} proved \eqref{addhermite} by computing the
corresponding Stieltjes tableau using induction.

For the Laguerre polynomials $L_n^{(\alpha)}(x)$ we have~\cite[p. 48]{Koe:Swa}:
\begin{align}
e^t\;
{}_0F_1(-; \alpha+1; -xt)=
\sum_{n=0}^\infty \frac{n!L_n^{(\alpha)}(x)}{(\alpha+1)_n}\frac{t^n}{n!}.
\end{align}
In \eqref{eq:addfor1F1} letting $B=\alpha+\beta+1$ and
substituting $t$ and $s$ by $t/\alpha$ and $s/\alpha$, respectively and, then
let $\alpha\to \infty$ we get
\bea
\begin{gathered}
{}_0F_1(-; B+1; -2(t+s))=\sum_{n= 0}^\infty
\frac{(-1)^n(B)_n(4ts)^n}{(B)_{2n}(B+1)_n}\;
{}_0F_1(-; B+2n+1; -2t)\\
\hspace{3cm}\times {}_0F_1(-; B+2n+1; -2s).
\end{gathered}
\eea
Let
\begin{align}
Q_n(t;\alpha)=\frac{t^n}{n!}\,e^t \; {}_0F_1(-; \alpha+2n+1; -2xt), \qquad n\geq 0.
\end{align}
Then we have the following addition formula.
\begin{thm}
The functions $\{Q_j(t; \alpha) \}$ have the addition formula
\begin{align}\label{eq:addLaguerre}
Q_0(t+s;\alpha)=\sum_{n=0}^\infty \frac{n!(\alpha)_n(-4x^2)^n}{(\alpha)_{2n}(\alpha+1)_n}
Q_n(t;\alpha)Q_n(s;\alpha).
\end{align}
\end{thm}
As an immediate consequence we have
$\lambda_n=\frac{n(\alpha+n-1)(-4x^2)}{(\alpha+2n-1)(\alpha+2n-2)(\alpha+n)}$ and
\begin{align}
H_{i,i+n}={n+i\choose i}\frac{n!L_n^{(\alpha+2i)}(x)}{(\alpha+2i+1)_n}.
\end{align}
For the Meixner polynomials $M_n(x;\beta,c)$ we have
\begin{align}
e^t \; {}_1F_1\left(\left.\begin{array}{c}
  -x \\
  \beta                \end{array}\right|
  \left(\frac{1-c}{c}\right)t
\right)=\sum_{n=0}^\infty \frac{M_n(x;\beta,c)}{n!}t^n.
\end{align}
In \eqref{eq:addfor1F1} substituting $\alpha+1$, $\alpha+\beta+2$,
$t$ and $s$
by $-x$, $\beta$, $(c-1)t/2c$ and $(c-1)s/2c$, respectively, we obtain
\bea
\begin{gathered}
{}_1F_1\left(\left. \begin{array}{c}
-x\\
\beta
\end{array}\right| \frac{1-x}{c}(t+s)\right)
= \Sum \frac{n!(-x)_n(\beta+x)_n(\beta-1)_n}
{(\beta-1)_{2n}(\beta)_n} \left(\frac{1-c}{c}\right)^{2n}
(ts)^n
\\
\times   {}_1F_1\left(\left. \begin{array}{c}
n-x\\
\beta+2n
\end{array}\right| \frac{1-c}{c}t\right)\;
{}_1F_1\left(\left. \begin{array}{c}
n-x\\
\beta+2n
\end{array}\right| \frac{1-c}{c}s\right).
\end{gathered}
\eea
Therefore define
$$
Q_n(x;\beta,c)=\frac{t^n}{n!}e^{t}{}_1F_1\left(\left. \begin{array}{c}
n-x\\
\beta+2n
\end{array}\right| \frac{1-c}{c}t\right)
$$
we have the following addition formula.
\begin{thm}
The functions $\{Q_j(t; \beta,c) \}$ have the addition formula
\begin{align}
Q_0(t+s;\beta,c)=\Sum \frac{n!(-x)_n(\beta+x)_n(\beta-1)_n}
{(\beta-1)_{2n}(\beta)_n} \left(\frac{1-c}{c}\right)^{2n}
Q_n(t;\beta,c)Q_n(s;\beta,c).
\end{align}
\end{thm}

In the same way we derive
$$
\lambda_n=\frac{n(-x+n-1)(\beta+x+n-1)(\beta+n-2)}
{(\beta+2n-2)(\beta+2n-3)(\beta+n-1} \left(\frac{1-c}{c}\right)^{2},
$$
and
\begin{align}
H_{i,i+n}={i+n\choose i}M_n(x-i;\beta+2i,c).
\end{align}

The Meixner-Pollaczek polynomials $P_n^{(\lambda)}(x;\phi)$ have the generating function
\begin{align}
 e^t{}_1F_1\left(\left.
\begin{array}{c}
\lambda+ix\\
2\lambda
\end{array}\right| (e^{-2i\phi}-1)t\right)
=\sum_{n=0}^\infty\frac{P_n^{(\lambda)}(x;\phi)}{(2\lambda)_ne^{in\phi}}t^n.
\end{align}
In \eqref{eq:addfor1F1} substituting $\alpha+1$, $\beta+1$,
$t$ and $s$
by $\lambda-1+ix$, $\lambda-1-ix$, $(1-e^{-2i\phi})t/2$ and $(1-e^{-2i\phi})s/2$, respectively, and letting
\begin{align}
 Q_n^{(\lambda)}(x;\phi)=\frac{t^n}{n!}e^t
{}_1F_1\left(\left. \begin{array}{c}
\lambda+ix+n\\
2\lambda+2n
\end{array}\right| (e^{-2i\phi}-1)t\right),
\end{align}
we obtain the following addition formula corresponding to Meixner-Pollaczek polynomials.
\begin{thm} We have the addition formula
\begin{align}
Q_0^{(\lambda)}(t+s;\phi)=\Sum \frac{n!(\lambda+ix)_n(\lambda-ix)_n(2\lambda-1)_n}
{(2\lambda-1)_{2n}(2\lambda)_n} 4^{n}
Q_n^{(\lambda)}(t;\phi)Q_n^{(\lambda)}(s;\phi).
\end{align}
\end{thm}

\section{$q$-Ultraspherical and Askey--Wilson Polynomials}
The continuous $q$-ultraspherical polynomials have the weight function
$$
w(x; \beta) =\frac{1}{2\pi}\frac{(e^{2i\theta},e^{-2i\theta} )}
{(\beta e^{2i\theta},\beta e^{-2i\theta})}\frac{(\beta^2,q)_\infty}
{(\beta,\beta q)_\infty}\frac{1}{\sqrt{1-x^2}},\quad x=\cos \theta \in (-1,1),
$$
and have the property
$$
C_n(x,\beta|q)=\frac{2^n(\beta)_n}{(q)_n}P_n(x),   \qquad
\lambda_j=\frac{(1-q^j)(1-\beta^2 q^{j-1})}{4(1-\beta q^{j-1})(1-\beta q^j)}.
$$
In view of  \eqref{eqUnasultra} and \eqref{eqPWUltra} we see that
\begin{align}
\label{eqPWExpinCheby}
e^{xy}
=\frac{2}{y}\sum_{n=0}^\infty (n+1)I_{n+1}(y)U_n(x).
\end{align}
Now the special case $\gamma = q$ of the connection relation \eqref{eq:connection} gives the following expansion:
\begin{align}
U_n(x)=\sum_{k=0}^{[n/2]}\frac{\beta^k(q/\beta;q)_k(q;q)_{n-k}}{(q;q)_k(q\beta;q)_{n-k}}
\frac{1-\beta q^{n-2k}}{1-\beta}C_{n-2k}(x;\beta\mid q).
\end{align}
Using the above expansion and the orthogonality of the $q$-ultraspherical polynomials
we have
\begin{align*}
Q_j(t)&=\frac{2}{t}\frac{1}{\lambda_1\cdots \lambda_j}\sum_{n=0}^\infty
\int_{-1}^1w(x)(n+1)I_{n+1}(t)U_n(x)\; P_j(x)dx\\
&=\frac{2^{j+1}\, (\beta)_j}{(q)_j}\frac{1}{t}
\sum_{k=0}^\infty (j+2k+1)I_{j+2k+1}(t)
\frac{\beta^k(q/\beta;q)_k(q;q)_{j+k}}
{(q;q)_k(q\beta;q)_{j+k}}\frac{1-\beta q^j}{1-\beta}\\
&=\frac{2^{j+1}}{t}
\sum_{k=0}^\infty (j+2k+1)I_{j+2k+1}(t)
\frac{\beta^k(q/\beta;q)_k(q^{j+1};q)_{k}}
{(q;q)_k(\beta q^{j+1};q)_{k}}.
\end{align*}

To denote the explicit dependence on $q$ and $\beta$ we set
\bea
\label{eqQiforqultra}
Q_j(t; \beta, q) = \frac{2^{j+1}}{t}
\sum_{k=0}^\infty (j+2k+1)I_{j+2k+1}(t)
\frac{\beta^k(q/\beta;q)_k(q^{j+1};q)_{k}}
{(q;q)_k(\beta q^{j+1};q)_{k}}.
\eea
Thus we proved that
\begin{thm}
The functions $\{Q_j(t; \beta, q) \}$ have the addition formula
\bea
\label{eqformthm5.1}
Q_0(s+t; \beta, q) = \Sum \frac{(q;q)_n(\beta^2;q)_n}{4^n(\beta;q)_n (q\beta;q)_n}
Q_n(s; \beta, q) Q_n(t; \beta, q).
\eea
\end{thm}
The special case $\beta \to 0$ is worth recording. Indeed if
\bea
F_n(t;q) :=  \frac{2^{n+1}}{t}
\sum_{k=0}^\infty (n+2k+1)I_{n+2k+1}(t) (-1)^n q^{\binom{k+1}{2}}\qbi{n+k}{k}{q},
\eea
then we have established  the curious result
\bea
F_0(s+t;  q) = \Sum \frac{(q;q)_n}{4^n}
F_n(s;  q) F_n(t; q).
\eea

Another interesting case is to let $\beta=q^\nu$ then let $q\to 1$.  This should reduce
\eqref{eqQiforqultra}  to \eqref{eqQiforultra} since $\lim_{q \to 1}C_n(x; q^\nu|q) = C_n^\nu(x)$. Surprisingly the $q \to 1$ limit of \eqref{eqQiforultra}, after setting $\beta = q^\nu$ is
\begin{align*}
Q_j(t)= 2^j\sum_{k=0}^\infty
I_{j+2k}(t)\frac{({-\nu})_k(j)_k}
{k!(j+1+\nu)_k}\frac{j+2k}{j}.
\end{align*}
Equating the above limit and the $Q_j$ as in  \eqref{eqQiforqultra} leads to the following known identity involving Bessel functions
\bea
\label{eqBesselidentity}
(z/2)^{\mu-\nu}J_\nu(z)=\sum_{n=0}^\infty\frac{\Gamma(\mu+n)\Gamma(\nu+1-\mu)(\mu+2n)}
{n!\,\Gamma(\mu+1-\mu-n)\Gamma(\nu+n+1)}J_{\mu+2n}(z),
\eea
see \cite[(7.15.2)]{Erd:Mag:Obe:Tri}. It is also worth mentioning that
\eqref{eqBesselidentity} is equivalent to a theorem of Bailey evaluating the sum of a
well-poised ${}_4F_3$ with argument $-1$, \cite[(4.5.4)]{Erd:Mag:Obe:Tri-1} .

Next we consider the Askey--Wilson polynomials whose weight function is
\bea
\begin{gathered}
W(x; a_1, a_2, a_3, a_4|q) = \frac{(e^{2i\theta}, e^{-2i\theta};q)_\infty}
{\prod_{j=1}^4(a_j e^{i\theta}, a_je^{-i\theta};q)_\infty}\; \frac{1}{\sqrt{1-x^2}}, \\
\times \frac{(q;q)_\infty \prod_{1\le j < k \le 4}(a_ja_k;q)_\infty} {2\pi \; (a_1a_2a_3a_4;q)_\infty},
\quad x = \cos \theta.
\end{gathered}
\eea
The \awp have the basic hypergeometric function representation
\begin{equation}
\label{eqAsWip}
\begin{split}
p_n(x; a_1, a_2, a_3, a_4\,|\, q) &=a_1^{-n}\(a_1a_2, a_1a_3, a_1a_4;q\)_n \\
&\quad\times {}_{4}\phi_3\(\left. \begin{matrix}
q^{-n}, a_1a_2a_3a_4q^{n-1}, a_1e^{i\theta}, a_1e^{-i\theta} \\
a_1a_2,\; a_1a_3, \;a_1a_4
\end{matrix}\, \right|q,q\).
\end{split}
\end{equation}
One very special case of their connection coefficients formula is \cite[(6.4)--(6.5)]{Ask:Wil}
\bea
p_n(x; \alpha, a_2, a_3, a_4\,|\, q) = \sum_{k=0}^nc_{k,n}p_n(x; a, a_2, a_3, a_4\,|\, q)
\eea
where
\bea
c_{k,n} = \frac{a^{n-k}(q;q)_n(\alpha a_2 a_3 a_4q^{n-1};  q)_k (\alpha/a;q)_{n-k}}
{(q, a a_2 a_3 a_4 q^{k-1} ; q)_k(q, a a_2 a_3 a_4q^{2k};q)_{n-k}}\; \prod_{2\le j <m \le 4}(a_ja_m q^k;q)_{n-k}.
\eea
Moreover
\bea
\begin{gathered}
U_n(x) = \frac{1}{(q^{n+2};q)_n}p_n(x; \sqrt{q}, q, -\sqrt{q}, -q|q), \quad p_n(x; a, b, c, d|q) = 2^n (abcdq^{n-1};q)_nP_n(x) \\
\lambda_n = \frac{(1-q^n)(1-a_1a_2a_3a_4 q^{n-2})\prod_{1 \le j < k \le 4}(1- a_j a_k q^{n-1})}
{4(a_1a_2a_3a_4 q^{2n-3}, a_1a_2a_3a_4 q^{2n-2};q)_2}.
\nonumber
\end{gathered}
\eea

Applying \eqref{eq:func} and the plane wave expansion \eqref{eqPWExpinCheby}  we find that the $Q_m$'s are given by
\bea
\begin{gathered}
Q_m(t) = \frac{1}{\lambda_1\lambda_2 \cdots \lambda_m} {\cal L} (e^{xt}P_m)
= \frac{2/t}{\lambda_1\lambda_2 \cdots \lambda_m} \Sum (n+1)I_{n+1}(t) {\cal L} ( U_n(x) P_m)\\
= (2/t) \sum_{n=m}^\infty  (n+1)I_{n+1}(t)
\frac{(aq^{m+1};q)_m}{(q^{n+2};q)_n} 2^n \; c_{m,n}\\
=\ (2/t) \sum_{n=m}^\infty 2^n  (n+1)I_{n+1}(t)
\frac{(aq^{m+1};q)_m}{(q^{n+2};q)_n}
 \frac{a^{n-m}(q;q)_n(q^{n+2}; q)_m (q/a;q)_{n-m}}
{(q, a q^{m+1} ;q)_m(q, aq^{2m+2};q)_{n-m}}\;
\\
\qquad \qquad \times (-q^{m+1}, q^{m+3/2}, -q^{m+3/2} ;q)_{n-m}.
\end{gathered}
\nonumber
\eea
After some simplification we arrive at
\bea
\label{eqQmforAW}
\begin{gathered}
Q_m(t) = \frac{2^{m+1}}{t }  \Sum
 2^n a^n (n+m+1)       \qquad \qquad \\
\qquad \qquad \times
 \frac{(q^{m+1}, q/a, -q^{m+1}     ;q)_n(q^{2m+3};q^2)_n}
 {(q, aq^{2m+2}, q^{n+m+2} ;q)_n} \; I_{n+m+1}(t).
\end{gathered}
\eea
\begin{thm}
The functions $\{Q_m(x)\}$ defined in \eqref{eqQmforAW} satisfy the addition theorem
\bea
\label{eqAdditionforQAW}
Q_0(s+t) = \Sum  \frac{(q^2, a^2q   ;q)_n}
{4^n(aq,  aq^2;q)_{2n}}\; Q_n(t)Q_n(s).
\eea
\end{thm}


\section{Ultraspherical Polynomials as Moments}
One of the generating functions reads
\begin{align}
Q(t)=\sum_{n=0}^\infty
\frac{C_{n}^{(\nu)}(x)}{(2\nu)_{n}}t^n=e^{xt}
{}_{0}F_{1}\left(
\begin{array}{c}
-\\
\nu+\frac{1}{2}
\end{array};\frac{(x^2-1)t^2}{4}
\right).
\end{align}
Let $\cos\phi=-1$, then $w=z+Z$ and $C_{n}^\nu(-1)=(-1)^n\frac{(2\nu)_{n}}{n!}$.
It follows from \eqref{eq:besseladd} that
\begin{align}
\begin{gathered}
{}_{0}F_{1}\left(
\begin{array}{c}
-\\
\nu+1/2
\end{array};\frac{-(z+Z)^2}{4}
\right)=
\sum_{n=0}^\infty
\frac{(n+\nu-1/2)(-1)^n(2\nu-1)_{n}}{n!(\nu+1/2)_{n}(\nu-1/2)_{n+1}}\\
\hfill\times \left(\frac{zZ}{4}\right)^n
{}_{0}F_{1}\left(
\begin{array}{c}
-\\
\nu+1/2+n
\end{array};\frac{-z^2}{4}
\right) {}_{0}F_{1}\left(
\begin{array}{c}
-\\
\nu+1/2+n
\end{array};\frac{-Z^2}{4}
\right).
\end{gathered}
\end{align}
Therefore, let $z=t\sqrt{1-x^2}$ and $Z=t\sqrt{1-x^2}$ we obtain
\begin{align}\label{eq:qultradd}
Q(t+s)=\sum_{n=0}^\infty
\frac{(n+\nu-1/2)(-1)^n(2\nu-1)_{n}}{(\nu+1/2)_{n}(\nu-1/2)_{n+1}}
\frac{(1-x^2)^nn!}{4^{n}}
Q_{n}(t)Q_{n}(s),
\end{align}
where
$$
Q_{n}(t)=\frac{t^n}{n!}e^{tx}{}_{0}F_{1}\left(
\begin{array}{c}
-\\
\nu+1/2+n
\end{array};\frac{(x^2-1)t^2}{4}
\right)=\frac{t^n}{n!}+(n+1)x\frac{t^{n+1}}{(n+1)!}+\cdots.
$$

Extracting the coefficients of $t^{m}s^{n}$ in \eqref{eq:qultradd} we get
\bea
\begin{gathered}
\frac{(m+n)!C_{m+n}^{\nu}(x)}
{(2\nu)_{m+n}m!n!}=\sum_{k=0}^\infty
\frac{(k+\nu-1/2)(-1)^k(2\nu-1)_{k}}
{k!(\nu+1/2)_{k}(\nu-1/2)_{k+1}}\frac{(1-x^2)^k}{4^{k}}
\\
\hspace{2cm}\times \frac{C_{m-k}^{\nu+k}(x)
C_{n-k}^{\nu+k}(x)}{(2\nu)_{m-k}(2\nu)_{n-k}}.
\end{gathered}
\eea
Using the relation
$$
\lim_{\alpha\mapsto \infty}\alpha^{-n/2}C_{n}^{\alpha+1/2}(x/\sqrt{\alpha})=\frac{H_{n}(x)}{n!},
$$
we derive
\begin{align}\label{eq:her}
\frac{H_{m+n}(x)}{m!n!}=\sum_{k= 0}^{m\wedge n} \frac{(-2)^k}{k!}\frac{H_{m-k}(x)}{(m-k)!}
\frac{H_{n-k}(x)}{(n-k)!}.
\end{align}
An immediate consequence of \eqref{eq:qultradd} is the following formula for the Hankel
determinant evaluation~\cite[Corollary~3]{Wim}.
\begin{cor}
We have
\begin{align}
\det\left(\frac{(i+j)!}{(2\nu)_{i+j}}C_{i+j}^{\nu}(x)\right)_{0\leq i,j\leq n}=
\frac{(x^2-1)^{n(n+1)/2}}{2^{n^2}}\prod_{r=1}^n\frac{r!(2\nu)_{r-1}}
{(\nu+1/2)_{r-1}(\nu+1/2)_r},
\end{align}
and more generally, for $n\geq 0$, the entries of the Stieltjes tableau~\eqref{stieltjes} are
\begin{align}
H_{i,i+n}=\sum_{k=0}^{\lfloor n/2\rfloor}\frac{(n+i)!\,x^{n-2k}(x^2-1)^k\,t^{2k}}
{i!k!(n-2k)!(\nu+1/2+i)_k 4^k}.
\end{align}
\end{cor}

\section{A variation of the Stieltjes-Rogers addition formula}
Let $\{P_n(x)\}$ satisfy \eqref{eq:threeterm} with moment sequence $\{\mu_n\}$
 and $\bar P_n(x)=a^{-n}P_n(ax+b)$
 ($a\neq 0$).
 Then it is well-known~\cite[p.25]{Chi} that
 $\{\bar P_n(x)\}$ is
an OPS with respect to the moments given by
\bea\label{eq:shiftmoment}
\bar \mu_n=a^{-n}\sum_{k=0}^n{n\choose k}(-b)^{n-k}\mu_k,
\eea
and satisfy
\begin{align}
 \bar P_{n+1}(x)=\left(x-\frac{b_n-b}{a}\right)\bar P_n(x)-\frac{\lambda_n}{a^2}\bar P_{n-1}(x).
\end{align}

Let $Q_0(t)=\sum_{n= 0}^\infty \mu_n \frac{t^n}{n!}$. Then it is easy to see that
\bea
\bar Q_0(t)=\sum_{n= 0} ^\infty \bar \mu_n \frac{t^n}{n!}=e^{-bt/a} Q_0(t/a).
\eea
The following  variation of the Stieltjes-Rogers addition formula \eqref{eqaddthm} is sometimes very useful.
\begin{thm}
The addition formula for the moment sequence \eqref{eq:shiftmoment} is
\begin{align}
\bar Q_0(s+t)=\sum_{n=0}^\infty \lambda_1\ldots \lambda_n a^{-2n}\bar Q_n(s)\bar Q_n(t),
\end{align}
where $\bar Q_n(t)=e^{-bt/a} Q_n(t/a)$. The corresponding entries in \eqref{stieltjes} are
\begin{align}
 \bar H_{j,j+n}=\sum_{k=0}^{n+j}{n+j\choose k}(-b)^ka^{-n-j}H_{j,n+j-k}.
\end{align}
\end{thm}
In particular we have $\bar H_{n,n}=(-b/a)^nH_{n,n}$, i.e.,
\begin{align}
\det(\bar \mu_{i+j})_{0\leq i,j\leq n}=\left(\frac{-b}{a}\right)^n\det( \mu_{i+j})_{0\leq i,j\leq n}.
\end{align}

For example, let $\mu_n=(\alpha+1)_n$ be the $n$th-moment of Laguerre polynomials
$\{L_n^{\alpha}(x)\}$ (see \eqref{laguerre}). If $a=b=1/x$ then
\begin{align}
\bar \mu_n=\sum_{k=0}^n(-1)^{n-k}{n\choose k}x^k (\alpha+1)_k
\end{align}
is a weighted derangement number.
The corresponding addition formula reads
\begin{align}
\bar Q_0(s+t)=\sum_{n=0}^\infty n! (\alpha+1)_n x^{2n}\bar Q_n(s)\bar Q_n(t),
\end{align}
The $\alpha=0$ case of the above formula was derived by Radoux~\cite{Rad} using induction.

\bigskip
{\bf Acknowledgments}:  We are grateful to Mizan Rahman for his simplification of the $Q_j$ functions of Sections  4 and 5 and to Tom Koornwinder for pointing out the work of Burchnal and Chaundy  \cite{Bur:Cha}. This work was done during the visit of the first author to Universit\'{e} Lyon 1 and he gratefully acknowledges the hospitality and financial support of the Institute Camille Jordan.

\end{document}